\input amstex
\magnification 1100
\documentstyle{amsppt}
\NoBlackBoxes

\def\Z{{\Bbb Z}}
\def\R{{\Bbb R}}
\def\C{{\Bbb C}}
\def\P{{\Bbb P}}
\def\calP{{\Cal P}}
\def\calK{{\Cal K}}
\def\t {{\frak t}}
\def\A {{\frak A}}
\def\diag{\operatorname{diag}}
\def\conv{\operatorname{conv}}

\def\refAW        {1}
\def\refBelkale   {2}
\def\refBertram   {3}
\def\refBCF       {4}
\def\refBuch      {5}
\def\refDK	  {6}	% Danilov-Koshevoi
\def\refFulton	  {7}
\def\refKlyachko  {8}
\def\refKT	  {9}
\def\refKTsurv    {10}
\def\refKTW	  {11}	% Knutson-Tao-Woodward
\def\refJKTR	  {12}
\def\refOO	  {13}  % Simplex-Method
\def\refTao       {14}

\def\eqLU	{2}

%%%%%%%%%%%%%%%%%%%%%%%%%%%%%%%%%%%%%%%%%%%%%%%%%%%%%%%%%%%%%%%%%%%%%%%%%%
%%%%%%%%%%%%%%%%%%%%%%%%%%%%%%%%%%%%%%%%%%%%%%%%%%%%%%%%%%%%%%%%%%%%%%%%%%
%%%%%%%%%%%%%%%%%%%%%%%%%%%%%%%%%%%%%%%%%%%%%%%%%%%%%%%%%%%%%%%%%%%%%%%%%%
%%%%%%%%%%%%%%%%%%%%%%%%%%%%%%%%%%%%%%%%%%%%%%%%%%%%%%%%%%%%%%%%%%%%%%%%%%
%%%%%%%%%%%%%%%%%%%%%%%%%%%%%%%%%%%%%%%%%%%%%%%%%%%%%%%%%%%%%%%%%%%%%%%%%%
%%%%%%%%%%%%%%%%%%%%%%%%%%%%%%%%%%%%%%%%%%%%%%%%%%%%%%%%%%%%%%%%%%%%%%%%%%

\leftheadtext { S.Yu.~Orevkov, Yu.P.~Orevkov }
\rightheadtext{ Agnihotri-Woodward-Belkale polytope and Klyachko cones }

\topmatter
\title
	Agnihotri-Woodward-Belkale polytope and the intersection
	of Klyachko cones
\endtitle
\author
	S.Yu.~Orevkov, $\boxed{\text{Yu.P.~Orevkov}}$
\endauthor
\address
	Laboratoire Emile Picard, UFR MIG, Univ. Paul Sabatier, Toulouse, France
	and Steklov Math. Inst., Moscow Russia
\endaddress
\address
	Faculty of Economy, Moscow State Univ., Moscow, Russia
\endaddress
\abstract
Agnihotri-Woodward-Belkale polytope $\Delta$ (resp. Klyachko cone $\calK$) is the 
set of solutions of the multiplicative (resp. additive) Horn's problem,
i.e., the set of triples of spectra of special unitary (resp. traceless
Hermitian) $n\times n$ matrices satisfying $AB=C$ 
(resp. $A+B=C$). $\calK$ is the tangent cone of $\Delta$ at the origin.
The group $G=\Z_n\oplus\Z_n$ acts naturally on $\Delta$.

In this note, we report on a computer calculation which shows that
$\Delta$ coincides with the intersection of $g\calK$, $g\in G$,
for $n\le14$ but does not coincide for $n=15$.

Our motivation was an attempt to understand how to solve the multiplicative 
Horn problem in practice for given conjugacy classes in $SU(n)$.
\endabstract

\endtopmatter

%%%%%%%%%%%%%%%%%%%%%%%%%%%%%%%%%%%%%%%%%%%%%%%%%%%%%%%%%%%%%%%%%%%%%%%%%%
%%%%%%%%%%%%%%%%%%%%%%%%%%%%%%%%%%%%%%%%%%%%%%%%%%%%%%%%%%%%%%%%%%%%%%%%%%
%%%%%%%%%%%%%%%%%%%%%%%%%%%%%%%%%%%%%%%%%%%%%%%%%%%%%%%%%%%%%%%%%%%%%%%%%%
%%%%%%%%%%%%%%%%%%%%%%%%%%%%%%%%%%%%%%%%%%%%%%%%%%%%%%%%%%%%%%%%%%%%%%%%%%
%%%%%%%%%%%%%%%%%%%%%%%%%%%%%%%%%%%%%%%%%%%%%%%%%%%%%%%%%%%%%%%%%%%%%%%%%%
%%%%%%%%%%%%%%%%%%%%%%%%%%%%%%%%%%%%%%%%%%%%%%%%%%%%%%%%%%%%%%%%%%%%%%%%%%

\document

\head Introduction
\endhead

For a special unitary matrix $A\in SU(n)$, let us denote by
$\lambda(A)$ its {\it spectrum} i.e. the vector
$(\lambda_1,\dots,\lambda_n)\in\R^n$ which is uniquely defined by the
conditions
$$
 \lambda_1 + \dots + \lambda_n = 0,
 \qquad
 \lambda_1\ge\lambda_2\ge\dots\ge\lambda_{n-1}\ge\lambda_n\ge\lambda_1-1,
        \eqno(1)
$$
and $\exp(2\pi i\lambda_1),\dots,\exp(2\pi i\lambda_n)$ 
are the eigenvalues of $A$. The mapping $\lambda:SU(n)\to\R^n$ identifies
conjugacy classes of $SU(n)$ with points of the 
$(n-1)$-simplex $\A$ given by (1). This simplex is called
{\it Weyl alcove} of $SU(n)$. It sits in the {\it Weyl chamber} $\t_+$
defined by (1) with the last inequality excluded.

Agnihotri-Woodward [\refAW] and Belkale [\refBelkale] obtained necessary and sufficient
conditions on three vectors $\alpha$, $\beta$, $\gamma$ for the existence of
matrices $A,B,C\in SU(n)$ such that $\alpha=\lambda(A)$, $\beta=\lambda(B)$,
$\gamma=\lambda(C)$, and $AB=C$.  They have shown that the
image of the mapping $\Lambda:SU(n)^2\to\A^3$, 
$(A,B)\mapsto(\lambda(A),\lambda(B),\lambda(AB))$, is a polytope
$\Delta\subset\A^3$  explicitely described
in terms of quantum Schubert calculus (see \S1).

This result is a generalization of Klyachko's solution [\refKlyachko]
of Horn's problem (see an excellent survey [\refFulton]).
Namely, Klyachko have described the set $\calK$ of
all triples of $n$-vectors realizable by spectra of traceless
Hermitian matrices $A$, $B$, and $A+B$. Of course, $\calK$ is the
tangent cone of $\Delta$ at the origin $O$ of $\R^{3n}$. 
Its facets are described in terms of the classical Schubert calculus.
We call $\Delta$ and $\calK$ the {\it Agnihotri-Woodward-Belkale
polytope} and the {\it Klyach\-ko cone} respectively.

How to apply these results in practice? Suppose, we have a concrete unitary matrix $C$
and we ask if it is realizable as a product of matrices of
given conjugacy classes (see [\refJKTR] for an example of applications).
Though there exists a finite set of inequalities to check,
the number of them grows exponentially and, for matrices, say, $30\times30$
there is no chance to generate all inequalities in a
reasonable time. This concerns both additive and multiplicative problems.
However, for the additive problem, the honeycomb description of $\calK$
due to Knutson and Tao [\refKT, \refKTsurv] allows to reduce the problem to the
existence of a solution of $3\binom{n}2$
simultaneous linear inequalities on $\binom{n-1}{2}$ variables (see also [\refDK]).
The same problem for $k$ matrices is reduced to the existence of a solution of
$3(k-2)\binom{n}2$ inequalities on $(k-2)\binom{n-1}2+(k-3)(n-1)$ variables.
For the multiplicative problem, such a reduction seems to be unknown (see [\refTao]).
This is why we study when the multiplicative problem can be
reduced to the additive problem.

Let $I$ be the identity $n\times n$ matrix and let $\omega=\exp(2\pi i/n)$.
Easy to see that $\A$ is the convex hull of $\lambda(Z)$ where
$Z=\{\omega^k I\}$ is the center of $SU(n)$.
Moreover, the action of $Z$ on $SU(n)$ by multiplication induces affine
linear actions of $Z$ on $\A$ and of $G=Z\times Z$ on $\Delta$.
This yields natural lower and upper estimates for
$\Delta$:
$$
	\conv(GO)\subset\Delta\subset\Delta_{\calK}
	\qquad\text{where}\quad
	\Delta_{\calK}=\bigcap_{g\in G} g\calK			\eqno(\eqLU)
$$
(``$\conv$'' stands for the convex hull; $GO=\{gO\,|\,g\in G\}=\Lambda(G)$).
Are these estimates sharp? For the lower estimate, this question was asked
and answered in [\refBelkale; Sect.~7]. The answer is: the equality
$\conv(GO)=\Delta$ holds for $n\le3$ and it fails for $n=4$.
For the upper estimate, we ask and answer this question in this note.
The answer is: the equality
$\Delta=\Delta_{\calK}$ holds for $n\le14$ and it fails for $n=15$.

Agnihotri and Woodward [\refAW; Sect.~8] discuss
the action of $G$ on $\Delta$ and its relation with a hidden symmetry
of Gromov-Witten invariants. In particular, they point out
that the action of $G$ allows to reduce degree $d$
invariants to degree zero invariants for small values of $n$
(which implies $\Delta=\Delta_{\calK}$ for those $n$), but,
for $n=10$, they give an example where the reduction is impossible.
In terms of the polyhedra $\Delta$ and $\Delta_{\Cal K}$, this means that
for $n=10$, among the inequalities which are used in [\refAW] for $\Delta$, there is
at least one which cannot be reduced by $G$ to a homogeneous inequality.
However, this inequality does not belong to a smaller system of inequalities
for $\Delta$ which is given by Belkale in [\refBelkale].
For the smaller system of inequalities, such examples occur only for $n\ge 15$.

%%%%%%%%%%%%%%%%%%%%%%%%%%%%%%%%%%%%%%%%%%%%%%%%%%%%%%%%%%%%%%%%%%%%%%%%%%
%%%%%%%%%%%%%%%%%%%%%%%%%%%%%%%%%%%%%%%%%%%%%%%%%%%%%%%%%%%%%%%%%%%%%%%%%%
%%%%%%%%%%%%%%%%%%%%%%%%%%%%%%%%%%%%%%%%%%%%%%%%%%%%%%%%%%%%%%%%%%%%%%%%%%
%%%%%%%%%%%%%%%%%%%%%%%%%%%%%%%%%%%%%%%%%%%%%%%%%%%%%%%%%%%%%%%%%%%%%%%%%%
%%%%%%%%%%%%%%%%%%%%%%%%%%%%%%%%%%%%%%%%%%%%%%%%%%%%%%%%%%%%%%%%%%%%%%%%%%
%%%%%%%%%%%%%%%%%%%%%%%%%%%%%%%%%%%%%%%%%%%%%%%%%%%%%%%%%%%%%%%%%%%%%%%%%%

\head 1. Description of $\Delta$ and $\calK$ in terms of (quantum)
	Schubert calculus
\endhead

Fix positive integers $r,k$ such that $r+k=n$.
Let $QH^*(G_r(\C^n))$ be the quantum cohomology ring of the
Grassmanian of $r$-planes in $\C^n$.
It is an algebra over $\Z[q]$ ($q$ is an indeterminate)
which is generated as a $\Z[q]$-module by the elements
$\{\sigma_a\}$ where $a$ runs over the set of partitions
$\calP_{r,k}=\{\,(a_1,\dots,a_r)\,|\,k\ge a_1\ge\dots\ge a_r\ge0\,\}$.
Let $N_{ab}^c(r,k)=\sum_{d=0}^\infty N_{ab}^{c,d}(r,k)\,q^d$ 
be the structure constants of this algebra, i.e.
$$
 \sigma_a\cdot\sigma_b = 
  \sum_{c\in\calP_{r,k}}\sum_{d=0}^\infty
  N_{ab}^{c,d}(r,k)\,\sigma_c q^d.
$$

The quantum multiplication is homogeneous if we set
$\deg q=n$, $\deg\sigma_a=|a|:=a_1+\dots+a_r$, i.e.
$N_{ab}^{c,d}(r,k)$ is nonzero only if $nd = |a|+|b|-|c|$.
If $d=0$ then $N_{ab}^{c,d}(r,k)$ 
coincides with the classical Littlewood-Richardson coefficient
$N_{ab}^c$ (in particular, it does not depend on $r$ and $k$).
An algorithm of computation of $N_{ab}^{c,d}(r,k)$ is given in [\refBCF].
It is implemented in [\refBuch].

\smallskip
Let
$\bar{\Cal I}=\{\, (r,k;a,b,c;d)\;\;|\;\; r+k=n,\;\; (a,b,c)\in\calP_{r,k}^3\,\}$.
For $t=(r,k;a,b,c;d)\in\Cal I$ we also denote $N_{ab}^{c,d}(r,k)$ by $N_t$.
Let $\Cal I=\{t\in\bar{\Cal I}\,|\, N_t=1\}$ and 
% Let $\Cal I$ be the subset of $\bar{\Cal I}$ 
% defined by $N_{ab}^{c,d}(r,k)=1$ and 
let $\Cal I_0$ be the subset of $\Cal I$ defined by $d=0$.

For $t=(r,k;a,b,c;d)\in\Cal I$, we
define $H_t=H_{ab}^{c,d}(r,k)$
as the half-space of $\R^{3n}$ given by the inequality 
$h_t(\alpha,\beta,\gamma)\ge 0$ where
$$
 h_t(\alpha,\beta,\gamma) = d + \sum_{i=1}^r \gamma_{k+i-c_i} 
 - \sum_{i=1}^r \alpha_{k+i-a_i}
 - \sum_{i=1}^r \beta_{k+i-b_i},
$$

As usually, we regard an elements $a=(a_1,\dots,a_r)$ of $\Cal P_{r,k}$
as a {\it Young diagram} inscribed in the rectangle $r\times k$ which has
$r$ rows (numbered from the top to the bottom) and $k$ columns
(numbered from the left to the right). The Young diagram of $a$ is the union
of $a_1$ leftmost squares of the first row, $a_2$ leftmost squares of the second
row etc. So, its area is $|a|$.

Let $\lambda:SU(n)\to\A$, $\Lambda:SU(n)^2\to\Delta$, $\t^+$,
and $\calK$ be as in Introduction.
The results of Klyachko and Agnihotri-Woodward-Belkale discussed in
Introduction state that 
$$
    \calK = \t_+^3\cap\Big(\bigcap_{t\in\Cal I_0} H_t\Big),
	\qquad
	\Delta = \A^3 \cap\Big(\bigcap_{t\in\Cal I} H_t\Big).
$$

% In this note we report on the results of computations which show
% that for $n\le 15$ each inequality (2)
% defines a facet of $\Delta_3$ when $N_{ab}^{c,d}(r,k)=1$, 
% i.e. the presented system of inequalities defining $\Delta_3$
% is minimal.

%%%%%%%%%%%%%%%%%%%%%%%%%%%%%%%%%%%%%%%%%%%%%%%%%%%%%%%%%%%%%%%%%%%%%%%%%%
%%%%%%%%%%%%%%%%%%%%%%%%%%%%%%%%%%%%%%%%%%%%%%%%%%%%%%%%%%%%%%%%%%%%%%%%%%
%%%%%%%%%%%%%%%%%%%%%%%%%%%%%%%%%%%%%%%%%%%%%%%%%%%%%%%%%%%%%%%%%%%%%%%%%%
%%%%%%%%%%%%%%%%%%%%%%%%%%%%%%%%%%%%%%%%%%%%%%%%%%%%%%%%%%%%%%%%%%%%%%%%%%
%%%%%%%%%%%%%%%%%%%%%%%%%%%%%%%%%%%%%%%%%%%%%%%%%%%%%%%%%%%%%%%%%%%%%%%%%%
%%%%%%%%%%%%%%%%%%%%%%%%%%%%%%%%%%%%%%%%%%%%%%%%%%%%%%%%%%%%%%%%%%%%%%%%%%

\head  2. The action of $G$
\endhead

Recall that $G=Z\times Z$ where $Z$ is the center of $SU(n)$.
Let $\Omega:\R^{n}\to\R^{n}$ be the affine linear mapping defined by
$$
\Omega(x_1,\dots,x_n) =
(x_2,x_3,\dots,x_n,\,x_1-1) + (1/n,\,\dots,\,1/n).
$$
It is easy to check that 
$\A = \bigcap_{j=0}^{n-1} \Omega^j(\t_+) = 
\conv\{\,\Omega^j(\vec0)\,|\,0\le j<n\}$.
The action of $G$ we spoke about in Introduction, is
$$
    g(\alpha,\beta,\gamma)=\big(\Omega^i(\alpha),\Omega^j(\beta),
    \Omega^{i+j}(\gamma)\big)\qquad\text{for}
    \quad g=(\omega^i,\omega^j).
$$
Obviously, $\Delta$ is invariant under $G$. Indeed,
$AB=C$ iff $(\omega^i A)(\omega^j B) = \omega^{i+j}C$.

It is shown in [\refAW; Sect.~8]
that the action of $G$ on $\Delta$ corresponds to the
action of $G$ on $\Cal I$ defined as follows.
Quantum Pierri formula (see [\refBertram]) implies that
$$
	\sigma_k\sigma_a = \cases
		\sigma_{k,a_1,\dots,a_{r-1}},	&a_r=0,\\
		\sigma_{a_1-1,\dots,a_r-1},	&a_r\ne0.
	\endcases
$$
Hence, for any $a,b,c\in \Cal P_{r,k}$ and $i,j\in\Z$, we have
$$
	\sigma_k^i\sigma_a = \sigma_{a'}q^{d'_a},
\quad
	\sigma_k^j\sigma_b = \sigma_{b'}q^{d'_b},
\quad
	\sigma_k^{i+j}\sigma_c = \sigma_{c'}q^{d'_c}
$$
for certain
$d'_a,d'_b,d'_c\ge0$ and $a',b',c'\in\calP_{r,k}$ uniquely
determined by $a,b,c$ and $i,j$. Thus, for 
$t=(r,k;a,b,c;d)$ and $t'=(r,k;a',b',c';d')$ where
$d'=d+d'_c-d'_a-d'_b$, we have
$ N_t = N_{t'} $.
In other words, the quantum Littlewood-Richardson coefficients are
symmetric under the action of $G$ on $\bar{\Cal I}$ defined by
$gt=t'$ for $g=(\omega^i,\omega^j)$
(this action is well defined because $\sigma_k^n = q^k$).
In particular, $\Cal I$ is invariant under the action of $G$.
It is shown in [\refAW] that the two actions of $G$ (on $\Cal I$ and
on $\Delta$) are coherent, i.e., $gH_t = H_{gt}$.

%%%%%%%%%%%%%%%%%%%%%%%%%%%%%%%%%%%%%%%%%%%%%%%%%%%%%%%%%%%%%%%%%%%%%%%%%%
%%%%%%%%%%%%%%%%%%%%%%%%%%%%%%%%%%%%%%%%%%%%%%%%%%%%%%%%%%%%%%%%%%%%%%%%%%
%%%%%%%%%%%%%%%%%%%%%%%%%%%%%%%%%%%%%%%%%%%%%%%%%%%%%%%%%%%%%%%%%%%%%%%%%%
%%%%%%%%%%%%%%%%%%%%%%%%%%%%%%%%%%%%%%%%%%%%%%%%%%%%%%%%%%%%%%%%%%%%%%%%%%
%%%%%%%%%%%%%%%%%%%%%%%%%%%%%%%%%%%%%%%%%%%%%%%%%%%%%%%%%%%%%%%%%%%%%%%%%%
%%%%%%%%%%%%%%%%%%%%%%%%%%%%%%%%%%%%%%%%%%%%%%%%%%%%%%%%%%%%%%%%%%%%%%%%%%

\head 3. Symmetries of $\Delta$ leaving $\calK$ invariant
\endhead

The the action of any non-trivial $g\in G$ is such that $g\calK\ne\calK$.
There are also evident symmetries which are common for $\Delta$ and $\calK$.

Since $AB=C \Leftrightarrow B^{-1}A^{-1}=C^{-1}$ and
$A+B=C\Leftrightarrow (-A)+(-B)=-C$, both
$\Delta$ and $\calK$ are invariant under the involution
$(\alpha,\beta,\gamma)\mapsto(\alpha^*,\beta^*,\gamma^*)$ where
$(\alpha_1,\dots,\alpha_n)^* = (-\alpha_n,\dots,-\alpha_1)$.
For $t\in\bar{\Cal I}$, we denote the image of $H_t$ under the mapping 
$(\alpha,\beta,\gamma)\mapsto(\alpha^*,\beta^*,\gamma^*)$ by $H^*_t$.
This action on the facets of $\Delta$ corresponds to
the following symmetry of (quantum) Littlewood-Richardson coefficients
provided by the isomorphism between the Grassmanians
$G_r(\C^n)$ and $G_k(\C^n)$.

Let $t\mapsto t^*$ be the involution $\bar{\Cal I}\to\bar{\Cal I}$ defined by
$$
	(r,k;\,a,b,c;\,d)^* = (k,r; a^*,b^*,c^*; d)
$$
where
$$
	(a_1,\dots,a_r)^*=(a_1^*,\dots,a_k^*),\qquad
	a_i^* = \max\{j\,|\, a_j\ge i \}
$$
(the Young diagrams of $a$ and $a^*$ are symmetric with respect to the
main diagonal).
For any $t\in\bar{\Cal I}$, we have $N_t = N_{t^*}$ and $H_t^* = H_{t^*}$.

For $a=(a_1,\dots,a_r)\in\P_{r,k}$, let $\bar a = (k-a_n,\dots,k-a_1)$.
Then $\bar{\Cal I}$ is invariant under the mappings
$t\mapsto(r,k;\,b,a,c;\,d)$ and $t\mapsto(r,k;\,b,\bar c,\bar a;\,d)$
(commutativity and the Poincar\'e duality).
These symmetries correspond to the evident symmetries of $\Delta$ provided by
$AB\sim BA$ and $AB=C \Leftrightarrow BC^{-1}=A^{-1}$.

Let $G_0$ be the group of linear transformations of $\R^{3n}$ generated
by all symmetries discussed in this section and let $\tilde G$ be the
group of affine transformations generated by $G$ and $G_0$. It is clear that
$|G_0|=12$ and $|\tilde G|=12n^2$.

%%%%%%%%%%%%%%%%%%%%%%%%%%%%%%%%%%%%%%%%%%%%%%%%%%%%%%%%%%%%%%%%%%%%%%%%%%
%%%%%%%%%%%%%%%%%%%%%%%%%%%%%%%%%%%%%%%%%%%%%%%%%%%%%%%%%%%%%%%%%%%%%%%%%%
%%%%%%%%%%%%%%%%%%%%%%%%%%%%%%%%%%%%%%%%%%%%%%%%%%%%%%%%%%%%%%%%%%%%%%%%%%
%%%%%%%%%%%%%%%%%%%%%%%%%%%%%%%%%%%%%%%%%%%%%%%%%%%%%%%%%%%%%%%%%%%%%%%%%%
%%%%%%%%%%%%%%%%%%%%%%%%%%%%%%%%%%%%%%%%%%%%%%%%%%%%%%%%%%%%%%%%%%%%%%%%%%
%%%%%%%%%%%%%%%%%%%%%%%%%%%%%%%%%%%%%%%%%%%%%%%%%%%%%%%%%%%%%%%%%%%%%%%%%%

\head 4. Statement of the results
\endhead

\proclaim{ Proposition 1 }
If $n\le14$, then $\Cal I = G\Cal I_0 = \bigcup_G g\Cal I_0$.
% $$
%  \{H_{ab}^{c,d}(r,k)\,|\,r+k=n,\, N_{ab}^{c,d}(r,k)=1\}
%  = \{gH_{ab}^{c,0}(r,k)\,|\,r+k=n,\, N_{ab}^c=1,\,g\in G\}.
%$$
In particular, $\Delta=\Delta_{\calK}$.
\endproclaim

\proclaim{ Proposition 2 } If $n=15$, then
$\Cal I=G\Cal I_0\cup \tilde G t_0$  and $t_0\not\in G\Cal I_0$ where
$$
	t_0=(6,9;\,663300,663300,666300;\,1).
$$
 In particular,
$\Delta=\Delta_{\calK}\cap
\big(\bigcap_{\tilde G} g H_{t_0} \big)$.

\medskip
Moreover, $\Delta\ne\Delta_{\calK}$, in particular 
$p=(\alpha,\beta,\gamma)\in(\Delta_{\calK}\setminus H_{t_0})
\subset (\Delta_{\calK}\setminus\Delta)$ for

\medskip
\centerline{
$\alpha=\beta={1\over17}\,(6,6,6,6,6,\;0,0,0,0,0,\;-6,-6,-6,-6,-6),
$}

\medskip
\centerline{
$
  \gamma={1\over17}\,(8,8,8,\;1,1,1,1,1,1,\;-3,-3,-3,-3,\;\,-9,-9).
$}
\medskip
The minimum of $h_{t_0}$ on $\Delta_{\calK}$
is attained at $p$.
\endproclaim

\proclaim{ Corollary 3 } If $n\le 15$, then the set of inequalities
$h_t\ge 0$, $t\in\Cal I$,
defining $\Delta$ is minimal, i.e., for any $t'\in\Cal I$, we have
$$
	\A^3\cap \Big(\bigcap_{t\in\Cal I\setminus\{t'\}} H_t\,\Big)\ne\Delta
$$
\endproclaim

\demo{ Proof } Follows from Propositions 1 and 2 combined with
Knutson-Tao-Woodward's result [\refKTW] on the minimality of the
system of inequalities $\{h_t\ge0\}_{t\in\Cal I_0}$ defining $\Cal K$.
\qed
\enddemo

\proclaim{ Proposition 4 } If $n=16$, then $\Cal I = G\Cal I_0 \cup \tilde G\Cal I'$ where
$$
\split
\Cal I'&=\{\,
	(6,10;\, 553300, 663300, 663300;\, 1),\\
	&(7,9;\, 5533000, 6633000, 6633000;\, 1),\,
	 (7,9;\, 5533300, 6633000, 6633300;\, 1),\\
	&(7,9;\, 5553000, 6443000, 6553000;\, 1),\,
	 (7,9;\, 6433300, 6633000, 6633300;\, 1),\\
	&(8,8;\, 44431000, 54441100, 54441100;\, 1),\,
	 (8,8;\, 44431000, 54442000, 54442000;\, 1),\\
	&(8,8;\, 44440000, 55441100, 55441100;\, 1),\,
	 (8,8;\, 44441000, 54441100, 55441100;\, 1),\\
	&(8,8;\, 44441100, 54431000, 54441100;\, 1),\,
	 (8,8;\, 44441100, 54441000, 55441100;\, 1),\\
	&(8,8;\, 44441100, 54441100, 55441110;\, 1)\,\}.
\endsplit
$$
\endproclaim

We do not know any human readable proof of Propositions 1, 2, and 4. They are
obtained on a computer. To check Proposition 1, we used the program
{\tt lrcalc} written by Buch [\refBuch]. Namely, using this program,
for each $n\le14$
we generated all elements $t\in\Cal I$  and, for each $t$, we checked
that its orbit $Gt$ (see \S2) contains an element with $d=0$.

Doing the same computations for $n=15$, we found the element $t_0$
whose orbit is disjoint from $\Cal I_0$. The only orbit
of $\tilde G$ disjoint from $\Cal I_0$ is $\tilde G t_0$.

To check that the inequality $h_{t_0}\ge 0$ is independent of the others,
we minimized (using [\refOO]) $h_{t_0}$ under
the constraints (1) and $h_t\circ g\ge0$, $t\in\Cal I_0$, $g\in G$.
These are 3\,135\,129\,030 constraints on 30 variables
(since $a=b$ in $t_0$, we may set $\alpha=\beta$).
The capacity of available computers was not enough for a straight forward
solution of this problem. In rest of this section we explain
the trick we used to reduce the
number of constraints up to 148\,295.

\medskip
Consider the following $15\times 15$ diagonal
matrices:
$$
\split
A&=\diag(1,1,1,1,1,\,-1,\dots,-1),\\
B_1&=\diag(1,1,1,\,-1,\dots,-1,\, 1,1),\\
B_2&=\diag(1,1,\, -1,\dots,-1,\, 1,1,1).
\endsplit
$$
Then 
\medskip\noindent\centerline{%
$\lambda(A)=\lambda(B_1)=\lambda(B_2)=
{1\over2}(-1,-1,-1,-1,-1,\,0,0,0,0,0,\,1,1,1,1,1)$,}
\medskip\noindent\centerline{%
$\lambda(AB_1) =
{1\over2}(-1,-1,\,0,\dots,0,\,1,1),$ and
$\lambda(AB_2)={1\over2}(-1,-1,-1,\,0,\dots,0,\,1,1,1)$.}
\medskip
\noindent
Let $p_i = \Lambda(A,B_i)$, % = (\lambda(A),\lambda(B_i),\lambda(AB_i))$,
$i=1,2$. By definition, we have $p_1,p_2\in\Delta$.
We have also $h_{t_0}(p_1) = h_{t_0}(p_2)=0$.
%Thus, the points $p_1$ and $p_2$ belong to the face $F_0$
%of $\Delta$ on which the minimum of $h_{t_0}$ is attained. 
%This face might have any dimension but we want to show that it is a facet
%(a hyperface) of $\Delta$, i.e., $\dim F_0 = \dim\Delta-1$.
Let 
$$
	\Delta' = \cap_{t\in\Cal I'} H_t,
	\qquad
	\Cal I' = \{t\in\Cal I\,|\, t\ne t_0,\; h_t(p_1)=h_t(p_2)=0\}
$$

\proclaim{ Lemma }
 If $\min_{\Delta'}h_{t_0} < 0$, then $\Delta\ne\Delta_{\calK}$.
\endproclaim

\demo{Proof}
Suppose that  $\min_{\Delta'}h_{t_0} < 0$. Let $p'\in\Delta'$ be
such that $h_{t_0}(p')<0$. Set $p_0=(p_1+p_2)/2$. Let 
$\Cal I'' = \Cal I\setminus(\Cal I'\cup\{t_0\})$ and
$Q=[p',p_0]\cap\Big(\bigcap_{t\in\Cal I''}\{h_t=0\}\Big)$.
For any $t\in\Cal I''$, the values of $h_t$ at $p_1$ and $p_2$ are
non-negative and at least one of them is positive, hence $h_t(p_0)>0$.
Therefore, $p_0\not\in Q$. Let $q$ be the point of $Q$ closest to $p_0$.
Then $h_t(q)\ge 0$ for any $t\in\Cal I''$. On the other hand, $h_t(q)\ge 0$
for any $t\in\Cal I'$ because $q\subset\Delta'$. Thus, 
$q\in\Delta_{\calK}\setminus\Delta$.
\qed\enddemo

Thus, to show that $\Delta\ne\Delta_{\Cal K}$, it is enough to find
the minimum of $h_{t_0}$ under the constraints only from $\Cal I'$ and this set
is more than 30\,000 times smaller than $\Cal I$.
The minimum is attained at the point $p=(\alpha,\beta,\gamma)$ presented
in Proposition 2 and we have $h_{t_0}(p)=-{1\over17}<0$. It follows from Lemma
that $\Delta\ne\Delta_{\calK}$. However, we cannot conclude that 
$\min_{\Delta_{\calK}}h_{t_0}=\min_\Delta h_{t_0}$. To prove this fact and to
ensure that floating point computations in [\refOO] did not affect the result,
we checked that $p\in\Delta_{\calK}$ (using {\tt lrcalc}, this takes few minutes).
Namely, we checked that
$h_t(p)<0$ only for $t=t_0$ and $h_t(p)=0$ only
in the following cases (up to swapping $a$ and $b$):
\roster
\item $r=3$, $a=b=(8,4,0)$, $c=(9,0,0)$;
\item $r=5$, $a=(7,7,3,0,0)$, $b=(7,7,4,4,0)$, $c=(8,8,6,1,1)$;
\item $r=7$, $a=(5,5,2,2,0,0,0)$, $b=(6,6,6,3,3,0,0)$, $c=(5,5,5,5,3,0,0)$;
\item $r=7$, $a=(5,5,3,3,3,0,0)$, $b=(5,5,3,3,3,0,0)$, $c=(5,5,5,5,3,0,0)$;
\item $r=13$, $a^*=b^*=(9,5)$, $c^*=(11,2)$;
\endroster
in all these cases we have $d=1$.

\Refs
\def\r{\ref}

\r\no\refAW
\by	S.~Agnihotri, C.~Woodward
\paper	Eigenvalues of products of unitary matrices and quantum Schubert
	calculus \jour Math. Research Letters \vol 5 \yr 1998 \pages 817--836
\endref

\r\no\refBelkale
\by	P.~Belkale
\paper	Local systems on $\Bbb P^1-S$ for $S$ a finite set
	\jour Compos. Math. \vol 129 \yr 2001 \pages 67--86
\endref

\r\no\refBertram
\by	A.~Bertram
\paper	Quantum Schubert calculus
	\jour Adv. Math., \vol 128 \yr 1997 \pages 289--305
\endref

\r\no\refBCF
\by	A.~Bertram, I.~Ciocan-Fontanine, W.~Fulton
\paper	Quantum multiplication of Schur polynomials
	\jour J. Algebra \vol 219 \yr 1999 \pages 728--746
\endref

\r\no\refBuch
\by	A.S.~Buch
\paper	{\tt lrcalc} -- Littlewood-Richardson Calculator
	\jour  Software available on \hbox to 2cm{}\newline
	{\tt http://home.imf.au.dk/$\widetilde{\;}$abuch/lrcalc} and on
	{\tt http://www-math.mit.edu/$\widetilde{\;}$abuch/lrcalc}
\endref

\r\no\refDK
\by	V.I.~Danilov, G.A.~Koshevoy
\paper	Discrete convexity and Hermitian matrices
\jour Tr. Mat. Inst. im. Steklova \vol 241 \yr 2003 \lang Russian \pages 68--89
\transl English transl.
\jour	Proc. Steklov Inst. Math. \vol 241 \yr 2003 \pages 58--78
\endref

\r\no\refFulton
\by	W.~Fulton
\paper	Eigenvalues, invariant factors, highest weights, and Schubert calculus
	\jour Bull. Amer. Math. Soc. (N.S.) \vol 37 \yr 2000 \pages 209--249
\endref

\r\no\refKlyachko
\by	A.A.~Klyachko
\paper	Stable bundles, representation theory and Hermitian operators
	\jour Selecta Math. (N.S.) \vol 4 \yr 1998 \pages 419--445
\endref

\r\no\refKT
\by	A.~Knutson, T.~Tao
\paper	The honeycomb model of $GL_n(\C)$ tensor products I:
	proof of the saturation conjecture
	\jour J. Amer. Math. Soc. \vol 12 \yr 1999 \pages 1055--1090
\endref

\r\no\refKTsurv
\by	A.~Knutson, T.~Tao
\paper	Honeycombs and sums of Hermitian matrices \jour Notices Amer. Math. Soc.
	\vol 48 \yr 2001 \issue 2 \pages 175--186
\endref

\r\no\refKTW
\by	A.~Knutson, T.~Tao, C.~Woodward
\paper	The honeycomb model of $GL_n(\C)$ tensor products II:
	puzzles determine facets of the Littlewood-Richardson cone
	\jour J. Amer. Math. Soc. \vol 17 \yr 2003 \pages 19--48
\endref

\r\no\refJKTR
\by	S.Yu.~Orevkov
\paper	Quasipositivity test via unitary representations of braid groups and its
	application to real algebraic curves
\jour	J. Knot Theory and Ramifications \vol 10 \yr 2001 \pages 1005--1023
\endref

\r\no\refOO
\by	S.Yu.~Orevkov, Yu.P.~Orevkov
\paper	Simplex-method {\rm(}linear programming{\rm)} \jour Software available on \newline
	{\tt http://picard.ups-tlse.fr/$\widetilde{\;}$orevkov/simplex}
\endref

\r\no\refTao
\by	T.~Tao
\paper	Open question: What is a quantum honeycomb? \hbox to 5cm{}\newline
	{\tt http://terrytao.wordpress.com/2007/04/19}
\endref

\endRefs

\enddocument